\documentclass[12pt]{article}
\usepackage{amsfonts}
\usepackage{amssymb}
\usepackage{amsthm}

\newtheorem{theorem}{Theorem}[section]

\newtheorem{lemma}{Lemma}[section]
\newtheorem*{remark}{Remark}

\newtheorem{definition}{Definition}[section]

\author{Yuri A. Turygin}
\title{A Borsuk-Ulam theorem for $(\mathbb Z_p)^k$-actions
on products of (mod $p$) homology spheres}
\date{}

\begin{document}
\maketitle

\begin{abstract}\scriptsize It is proved that for a product action of $(\mathbb
Z_p)^k$ on a product of (mod p) homology spheres
$N^{n_1}\times...\times N^{n_k}$, where all $n_i$'s are assumed to
be odd if $p$ is odd, and any continuous map $f\colon
N^{n_1}\times...\times N^{n_k}\to \mathbb R^m$ the set
$A(f)=\{x\in N^{n_1}\times...\times N^{n_k}|\ f(x)=f(gx)\ \forall
g\in(\mathbb Z_p)^k\}$ has dimension at least
$n_1+...+n_k-m(p^k-1)$, provided $n_i\ge mp^{i-1}(p-1)$ for all
$i\  (1\le i\le k)$. Moreover, if $n_i\ge mp^{k-1}(p-1)$ for all
$i(1\le i\le k)$ then the free action $\mu$ can be assumed
arbitrary.
\end{abstract}

\section{Introduction}
The famous Borsuk-Ulam theorem has been generalized by a number of
authors. One of the first and memorable generalizations is due to
C.T. Yang \cite{CTY2} and D.G. Bourgin \cite{B}:

\begin{theorem} Let $T$ be a fixed point free involution on a
sphere $S^n$ and let $f\colon S^n\to\mathbb R^m$ be a continuous
map into Euclidean space. Then the dimension of the coincidence
set $A(f)=\{x\in S^n|f(x)=f(Tx)\}$ is at least $n-m$.
\end{theorem}

The theorem above was generalized in works of several authors.
First, P.E. Conner and E.E. Floyd \cite{CF} proved the following
generalization:

\begin{theorem} Let $T$ be a differentiable involution on a sphere
$S^n$ and let $f\colon S^n\to M^m$ be a continuous map into a
differentiable manifold $M^m$ of dimension $m$. Assume that
$f_*\colon H_n(S^n;\mathbb Z_2)\to H_n(M;\mathbb Z_2)$ is trivial.
Then the dimension of the coincidence set $A(f)=\{x\in
S^n|f(x)=f(Tx)\}$ is at least $n-m$.
\end{theorem}

In their consequent works H. Munkholm \cite{M1} and M. Nakaoka
\cite{N} showed that the differentiability condition on the
involution $T$ can be dropped provided the target topological
manifold $M^m$ is assumed to be compact. Moreover, they
generalized the previous theorem to the case of free actions of a
cyclic group $\mathbb Z_p$ on (mod p) homology spheres. They
proved the following:

\begin{theorem} Let a cyclic group $\mathbb Z_p$ of a prime order
act freely on a (mod p) homology $n$-sphere $N^n$, and let
$f\colon N^n\to M^m$ be a continuous map into a compact
topological manifold $M^m$ of dimension $m$. If $p$ is odd also
assume that $M$ is orientable. Suppose that $f_*\colon
H_n(N;\mathbb Z_p)\to H_n(M;\mathbb Z_p)$ is trivial. Then the
dimension of the coincidence set $A(f)=\{x\in N|f(x)=f(gx)|\
\forall g\in\mathbb Z_p\}$ is at least $n-m(p-1)$.
\end{theorem}

One of the corollaries of the theorem above is the famous theorem
by J. Milnor \cite{Milnor} which asserts that every element of
order two in a group which acts freely on a sphere must be central
(see \cite{N} for details).

The purpose of this paper is to suggest another generalization of
the Borsuk-Ulam theorem, namely, to prove the next theorem.
Further and until the rest of the paper $p$ is always assumed to
be a prime number.

\begin{theorem}\label{T} Let $M:=N^{n_1}\times...\times N^{n_k}$ be
a product of (mod $p$) homology $n_i$-spheres and let $\mu\colon
(\mathbb Z_p)^k\circlearrowleft M$ be a product of free actions
$\mu_i\colon\mathbb Z_p\circlearrowleft N^{n_i}\  (1\le i\le k)$.
If $p$ is odd also assume that all $n_i$'s are odd. For a map
$f\colon M\to\mathbb R^m$ define a coincidence set $A(f):=\{x\in
M|f(x)=f(gx)\  \forall g\in (\mathbb Z_p)^k \}$. Then $$\dim
A(f)\ge \dim M - m(p^k-1)$$ provided $n_i\ge mp^{i-1}(p-1)$ for
all $i(1\le i\le k)$. Moreover, if one assumes $n_i\ge
mp^{k-1}(p-1)$ for all $i(1\le i\le k)$ then the free action $\mu$
can be assumed arbitrary.
\end{theorem}

\begin{remark} For $p=2$ and $m=1$ the theorem above was
implicitly proved by A.N. Dranishnikov in \cite{AD2}. In the case
$n_i\ge m(p^k-1)$ for all $i(1\le i\le k)$ the theorem above was
proved by V.V. Volovikov in \cite{Vol}.
\end{remark}

Let $G$ be a group and let $R$ be a commutative ring with one.
Then by $I_R(G)$ we denote the augmentation ideal of the group
ring $R[G]$, i.e. the kernel of the augmentation homomorphism
$R[G]\to R$.

The key ingredient in the proofs of the most Borsuk-Ulam type
theorems for maps into Euclidean spaces is the following basic
observation:

\begin{lemma} \label{Lemma}
Let $G\circlearrowleft M$ be a free action of a finite
group $G$ on a topological manifold $M$. For a continuous map
$f\colon M\to\mathbb R^m$ define a coincidence set $A(f):=\{ x\in
M|f(x)=f(gx)\  \forall g\in G\}$. Then $A(f)\neq\varnothing$ if
and only if the vector bundle $\xi\colon M\times_G I_{\mathbb
R^m}(G)\to M/G$ does not have a non-vanishing section.
\end{lemma}

\begin{proof} First, note that every continuous map $f\colon M\to\mathbb
R^m$ gives rise to a continuous section $\hat s(f)\colon M/G\to
M\times_G \mathbb R^m[G]$ of a vector bundle $\hat\xi\colon
M\times_G \mathbb R^m[G]\to M/G$ defined by a formula:
$$\hat s(f)(xG)=(x,\sum_{g\in G}f(xg^{-1})g)G.$$
Observe that $\hat\xi=\xi\oplus\varepsilon^m_{\mathbb R}$ where
$\varepsilon^m_{\mathbb R}$ is a trivial $m$-dimensional real
vector bundle. Therefore a projection $\pi\colon M\times_G \mathbb
R^m[G]\to M\times_G I_{\mathbb R^m}(G)$ is well defined. Now
define a continuous section $s(f)\colon M/G\to M\times_G
I_{\mathbb R^m}(G)$ of $\xi$ by a formula $s(f):=\pi\circ\hat
s(f)$. It is easy to see that $s(f)(xG)=0$ if and only if the
orbit of $x\in M$ is mapped by $f$ to a point.

Conversely, given a continuous section $s$ of $\xi$, it defines a
$G$-equivariant map $\bar s\colon M\to M\times \mathbb R^m[G]$
which is due to its equivariance must be of the form $\bar
s(x)=(x,\sum_{g\in G}f(xg^{-1})g)$ for some $f\colon M\to\mathbb
R^m$, and the lemma follows.
\end{proof}

Usually, to prove a Borsuk-Ulam type theorem for maps into
Euclidean spaces one shows that the Euler class of the vector
bundle $\xi\colon M\times_G I_{\mathbb R^m}(G)\to M/G$ in a
suitable cohomology theory is non-trivial. Then the dimension
restrictions on the coincidence set $A(f)$ follow (see the proof
of Theorem \ref{T}). For instance, the theorems from
\cite{M1,M,Roberts} were proved in this way. Unfortunately, when
one uses ordinary cohomology theory, Euler class of $\xi$ very
often turns out to be trivial (see \cite{M1}). This in fact is an
explanation of why all available results in the area are
restricted to the actions of so few groups. In this setting the
results of H. Munkholm from \cite{M} are especially interesting.
In that paper he proves a Borsuk-Ulam type theorem for $\mathbb
Z_{p^a}$-actions, $p$ is odd, on odd dimensional spheres using a
$\widetilde{KU}$-theory Euler class.

The proof the Theorem \ref{T} also will be based on the
non-triviality of the (mod $p$) Euler class of a corresponding
vector bundle. The next two sections will be devoted to the
calculation of Euler classes of relevant vector bundles.

\section{Calculation of $w_{2^k-1}(\eta)$}\label{mod2}

In this section assume that $G=(\mathbb Z_2)^k$. As usual $BG$
stands for the classifying space of $G$ and $EG$ stands for the
total space of the universal $G$-bundle. This section is devoted
to the calculation of the (mod $2$) Euler class of a vector bundle
$\eta\colon EG\times_G I_{\mathbb R}(G)\to BG$, i.e. its
Stiefel-Whitney class $w_{2^k-1}(\eta)$. These calculations are
then needed in the proof of Theorem \ref{T} in case $p=2$. Recall
that $H^*(BG;\mathbb Z_2)$ is a polynomial algebra $\mathbb
Z_2[x_1,...,x_k]$ on $1$-dimensional generators.

\begin{lemma}
\label{xi2} $w_{2^k-1}(\eta)=\prod_{q=1}^k\prod_{1\le
i_1<...<i_q\le k} (x_{i_1}+...+x_{i_q})$
\end{lemma}

\begin{proof} Let $\mathbb Z_2$ act on $\mathbb R$ by an obvious
involution. This involution induces on $\mathbb R$ a structure of
an $\mathbb R[\mathbb Z_2]$-module which we will denote by $V$.
Denote by $pr_i\colon BG\to \mathbb RP^\infty$ a projection on the
$i^{th}$ coordinate. Then by $\lambda_i$ we denote a
$1$-dimensional real vector bundle obtained from the following
diagram:

$$
\begin{array}{ccc}
E(\lambda_i)&\longrightarrow&S^\infty\times_{\mathbb Z_2} V\\
\lambda_i\downarrow&&\downarrow\\
BG&\stackrel{pr_i}{\longrightarrow}&\mathbb RP^\infty\\
\end{array}
$$
\\
Here $S^\infty$ stands for the infinite dimensional sphere. From
the construction of $\lambda_i$ it follows that
$w_1(\lambda_i)=x_i$.

Let $\eta_i$ be a vector bundle obtained from the following
diagram:

$$
\begin{array}{ccc}
E(\eta_i)&\longrightarrow&S^\infty\times_{\mathbb Z_2} \mathbb R[\mathbb Z_2]\\
\eta_i\downarrow&&\downarrow\\
BG&\stackrel{pr_i}{\longrightarrow}&\mathbb RP^\infty\\
\end{array}
$$
\\
From the equality $\mathbb R[\mathbb Z_2]=V\oplus
 (V\otimes_{\mathbb R[\mathbb Z_2]} V)=V\oplus V^2$ it follows that
$\eta_i=\lambda_i\oplus\lambda_i^2$ where
$\lambda_i^2=\lambda_i\otimes\lambda_i$ is a trivial
$1$-dimensional bundle. Recall an isomorphism of $\mathbb
R$-modules: $\mathbb R[G]\simeq\mathbb R[\mathbb
Z_2\oplus...\oplus\mathbb Z_2]\simeq\mathbb R[\mathbb
Z_2]\otimes_{\mathbb R}...\otimes_{\mathbb R}\mathbb R[\mathbb
Z_2].$ From this isomorphism it follows that
$\eta\oplus\varepsilon_{\mathbb
R}^1=\eta_1\otimes...\otimes\eta_k$. Therefore there exists the
following chain of isomorphisms of vector bundles:
$$\eta\oplus\varepsilon^1_{\mathbb R}\simeq \bigotimes_{i=1}^k
(\lambda_i\oplus\lambda_i^2)\simeq\bigoplus_{(\alpha_1,...,\alpha_k)\in
G}(\lambda_1^{\alpha_1}\otimes...\otimes\lambda_k^{\alpha_k}).$$

It is a well known that the first Stiefel-Whitney class of a
tensor product of $1$-dimensional real vector bundles equals to
the sum of the first Stiefel-Whitney classes of the multiplies.
Then by this fact and a formula of Whitney we get the following
chain of equalities:

$$w_{2^k-1}(\eta)=w_{2^k-1}(\eta\oplus\varepsilon^1_{\mathbb R})=
\prod_{(\alpha_1,...,\alpha_k)\neq 0}(\alpha_1 x_1+...+\alpha_k
x_k)=$$
$$=\prod_{q=1}^k\prod_{1\le i_1<...<i_q\le k} (x_{i_1}+...+x_{i_q}).$$
\end{proof}

\section{Euler class of $\eta_{\mathbb C}\colon EG\times_G I_{\mathbb C}(G)\to
BG$}\label{modp}

Through out this section assume that $p$ is a fixed odd prime and
that $G=(\mathbb Z_p)^k$. In this section we will calculate the
(mod $p$) Euler class of a complex vector bundle $\eta_{\mathbb
C}\colon EG\times_G I_{\mathbb C}(G)\to BG$ which equals to its
Chern class $c_{p^k-1}(\eta_{\mathbb C})$. These calculations are
then needed in the proof of Theorem \ref{T} in case of odd primes.
Recall that:
$$H^*(BG;\mathbb Z_p)=\Lambda_{\mathbb
Z_p}(y_1,...,y_k)\otimes\mathbb Z_p[x_1,...,x_k],$$ where
$\Lambda_{\mathbb Z_p}(y_1,...,y_k)$ is an exterior algebra on
$1$-dimensional generators and $\mathbb Z_p[x_1,...,x_k]$ is a
polynomial algebra on $2$-dimensional generators.

Chern classes of a regular representation of $G$, i.e. Chern
classes of the vector bundle $\eta_{\mathbb
C}\oplus\varepsilon_{\mathbb C}^1\colon EG\times_G \mathbb C[G]\to
BG$, were first computed by B.M. Mann and R.J. Milgram in
\cite{MM}. The lemma which is stated after the next definition is
essentially borrowed from their paper.

\begin{definition}
$L_k=\prod_{i=1}^k \prod_{\alpha_j\in \mathbb Z/p} (\alpha_1
x_1+...+\alpha_{i-1} x_{i-1}+x_i)$
\end{definition}

The polynomial defined above is called the $k^{th}$ Dickson's
polynomial (see \cite {MM} for more details).

\begin{lemma}
\label{euler} $e(\eta_{\mathbb C})=(-1)^k L_k^{p-1}$
\end{lemma}

\begin{proof}
The action of $\mathbb Z_p$ on $\mathbb C$ by rotations by
$\frac{2\pi}{p}$ induces on $\mathbb C$ a structure of a $\mathbb
C[\mathbb Z_p]$-module which we will denote by $L$. Let
$pr_i\colon BG\to B\mathbb Z_p$ be a projection on the $i^{th}$
coordinate. Then let $\lambda_i$ be a $1$-dimensional complex
vector bundle obtained from the following diagram:

$$
\begin{array}{ccc}
E(\lambda_i)&\longrightarrow&S^\infty\times_{\mathbb Z_p} L\\
\lambda_i\downarrow&&\downarrow\\
BG&\stackrel{pr_i}{\longrightarrow}&B\mathbb Z_p\\
\end{array}
$$
\\
It is not very difficult to show that $c_1(\lambda_i)=x_i$.

Let $\eta_i$ be a vector bundle obtained from the following
diagram:

$$
\begin{array}{ccc}
E(\eta_i)&\longrightarrow&S^\infty\times_{\mathbb Z_p} \mathbb C[\mathbb Z_p]\\
\eta_i\downarrow&&\downarrow\\
BG&\stackrel{pr_i}{\longrightarrow}&B\mathbb Z_p\\
\end{array}
$$
\\
It follows from the equality $\mathbb C[\mathbb
Z_p]=L\oplus...\oplus L^p$, where $L^j=L\otimes_{\mathbb C[\mathbb
Z_p]}...\otimes_{\mathbb C[\mathbb Z_p]} L$, that
$\eta_i=\lambda_i\oplus...\oplus\lambda_i^p$. Here
$\lambda_i^j=\lambda_i\otimes_{\mathbb C}...\otimes_{\mathbb
C}\lambda_i$. Also note that $\lambda_i^p$ is a trivial
$1$-dimensional complex bundle. Recall an isomorphism of $\mathbb
C$-modules: $\mathbb C[G]\simeq\mathbb C[\mathbb
Z_p\oplus...\oplus\mathbb Z_p]\simeq\mathbb C[\mathbb
Z_p]\otimes_{\mathbb C}...\otimes_{\mathbb C}\mathbb C[\mathbb
Z_p].$ From this isomorphism it follows that $\eta_{\mathbb
C}\oplus\varepsilon_{\mathbb C}^1=\eta_1\otimes...\otimes\eta_k$.
Therefore there exists the following chain of isomorphisms of
vector bundles:

$$\eta_{\mathbb C}\oplus\varepsilon^1_{\mathbb C}\simeq \bigotimes_{i=1}^k
(\lambda_i\oplus...\oplus\lambda_i^p)\simeq\bigoplus_{(\alpha_1,...,\alpha_k)\in
G}(\lambda_1^{\alpha_1}\otimes...\otimes\lambda_k^{\alpha_k}).$$

From a formula by Whitney and the fact that the first Chern class
of a tensor product of $1$-dimensional complex bundles equals to
the sum of the first Chern classes of the multiples, it follows
that:

$$c_{p^k-1}(\eta_{\mathbb C})=c_{p^k-1}(\eta_{\mathbb C}\oplus\varepsilon^1_{\mathbb
C})=\prod_{(\alpha_1,...,\alpha_k)\neq 0}(\alpha_1
x_1+...+\alpha_k x_k)=$$
$$=\prod_{i=1}^k [(p-1)!]^{p^{k-i}}
\prod_{(\alpha_1,...,\alpha_{i-1},1,0,...,0)}(\alpha_1
x_1+...+\alpha_{i-1} x_{i-1}+x_i)^{p-1}=$$
$$=[(p-1)!]^k L_k^{p-1}=(-1)^k L_k^{p-1}.$$
The last equality follows from a theorem of Wilson which states
that $(p-1)!\equiv (-1) (mod\ p)$. Thus $e(\eta_{\mathbb
C})=c_{p^k-1}(\eta_{\mathbb C})=(-1)^k L_k^{p-1}$.
\end{proof}

\section{Proof of Theorem \ref{T}}

In this section we will prove the main result of the paper Theorem
\ref{T}. Here assume that $p$ is any prime number and that
$G=(\mathbb Z_p)^k$.

\begin{proof}[Proof of Theorem \ref{T}]

Recall that $M=N^{n_1}\times...\times N^{n_k}$ is a product of
(mod $p$) homology $n_i$-spheres. We will begin the proof by
showing that under assumptions of the theorem the (mod $p$) Euler
class of $\xi_M\colon M\times_{\mu} I_{\mathbb R^m}(G)\to M/G$ is
non-trivial.

By universality property there exists the following commutative
diagram:

$$
\begin{array}{ccc}
M\times_{\mu} I_{\mathbb R^m}(G)&\longrightarrow&EG\times_G I_{\mathbb R^m}(G)\\
\xi_M\downarrow&&\downarrow\xi\\
M/G&\stackrel{\varphi}{\longrightarrow}&BG.\\
\end{array}
$$

{\sl Case $p$=2.} Let $\eta$ be a vector bundle from section
\ref{mod2}. Then from an isomorphism $I_{\mathbb R^m}=I_{\mathbb
R}\oplus...\oplus I_{\mathbb R}=mI_{\mathbb R}$ it follows that
$\xi=\eta\oplus...\oplus\eta=m\eta$. Thus
$e_2(\xi)=w_{2^k-1}(\xi)=w_{2^k-1}(\eta)^m$. By Lemma \ref{xi2} we
have $$w_{2^k-1}(\eta)=\prod_{q=1}^k\prod_{1\le i_1<...<i_q\le k}
(x_{i_1}+...+x_{i_q})=$$
$$=x_k^{2^{k-1}}\prod_{q=1}^{k-1}\prod_{1\le i_1<...<i_q\le k-1}
(x_{i_1}+...+x_{i_q})+R_k,$$ where $R_k$ contains monomials in
powers less than $2^{k-1}$. Therefore
$$e_2(\xi)=x_1^mx_2^{2m}\cdot...\cdot x_k^{2^{k-1}m}+Q_k,\eqno (1)$$ where
$Q_k$ does not contain monomials of the form
$x_1^mx_2^{2m}\cdot...\cdot x_k^{2^{k-1}m}$. It is easy to verify
that
$$H^*(M/G;\mathbb Z_2)=\mathbb
Z_2[x_1,...,x_k]/(x_1^{n_1+1},...,x_k^{n_k+1}),$$ and
$\varphi^*\colon H^*(BG;\mathbb Z_2)\to H^*(M/G;\mathbb Z_2)$ is
an epimorphism with $$Ker\
\varphi^*=(x_1^{n_1+1},...,x_k^{n_k+1}).$$ Thus from (1) and the
assumption $n_i\ge m2^{i-1}$ for all $i(1\le i\le k)$ it follows
that $e_2(\xi_M)=\varphi^*(e_2(\xi))\neq 0$.

{\sl Case $p>2$.} Let $\eta_{\mathbb C}$ be a vector bundle from
section \ref{modp}. Then from an isomorphism $I_{\mathbb
C^m}=I_{\mathbb C}\oplus...\oplus I_{\mathbb C}=mI_{\mathbb C}$ it
follows that $\mathbb C\xi=\eta\oplus...\oplus\eta=m\eta$, where
$\mathbb C\xi$ is a complexification of the vector bundle $\xi$.
We have the following chain of equalities:
$$e_p(\xi)^2=e_p(\mathbb C\xi)=c_{2^k-1}(\mathbb
C\xi)=c_{2^k-1}(\eta_{\mathbb C})^m.\eqno(2)$$ By Lemma
\ref{euler} we have
$$e_p(\eta_{\mathbb C})=(-1)^kL_k^{p-1}=(-1)^kL_{k-1}^{p-1}\left[ \prod_{\alpha_j\in\mathbb Z_p}
(\alpha_1x_1+...+\alpha_{k-1}x_{k-1}+x_k)\right]^{p-1}=$$
$$=(-1)^kx_k^{p^{k-1}(p-1)}L_{k-1}^{p-1}+R_k,$$ where $R_k$ contains
$x_k$ in powers less than $p^{k-1}(p-1)$. Thus
$$e_p(\mathbb C\xi)=
(-1)^{km}x_k^{mp^{k-1}(p-1)}L_{k-1}^{m(p-1)}+\hat R_k,$$ where
$\hat R_k$ contains $x_k$ in powers less than $mp^{k-1}(p-1)$.
Then by induction it follows that

$$e_p(\mathbb C\xi)=(-1)^{km}x_1^{m(p-1)}x_2^{mp(p-1)}\cdots x_k^{mp^{k-1}(p-1)}
+Q_k,$$ where $Q_k$ contains no monomials of the form
$$bx_1^{m(p-1)}x_2^{mp(p-1)}\cdots x_k^{mp^{k-1}(p-1)},\  b\neq 0,\  b\in\mathbb Z_p.$$
Therefore from the previous and (2) it follows that
$$e_p(\xi)=ax_1^{\frac{m(p-1)}{2}}x_2^{\frac{mp(p-1)}{2}}\cdots x_k^{\frac{mp^{k-1}(p-1)}{2}}
+\hat Q_k,\eqno (3)$$ where $a^2\equiv(-1)^q(mod\  p)$ for some
$q\ge 0$ and $\hat Q_k$ contains no monomials of the form
$$bx_1^{\frac{m(p-1)}{2}}x_2^{\frac{mp(p-1)}{2}}\cdots
x_k^{\frac{mp^{k-1}(p-1)}{2}},\  b\neq 0,\  b\in\mathbb Z_p.$$

It is not very difficult to see that $$H^*(M/G;\mathbb
Z_p)=\Lambda_{\mathbb Z_p}(y_1,...,y_k)\otimes_{\mathbb
Z_p}\mathbb
Z_2[x_1,...,x_k]/(x_1^{\frac{n_1+1}{2}},...,x_k^{\frac{n_k+1}{2}}),$$
where $\dim x_i=2$, and $\varphi^*\colon H^*(BG;\mathbb Z_p)\to
H^*(M/G;\mathbb Z_p)$ is an epimorphism with $$Ker\
\varphi^*=(x_1^{\frac{n_1+1}{2}},...,x_k^{\frac{n_k+1}{2}}).$$
Then from (3) and the assumption $n_i\ge mp^{i-1}(p-1)$ for all $i
(1\le i\le k)$ it follows that $e_p(\xi_M)=\varphi^*(e_p(\xi))\neq
0$.

Since $A(f)$ is closed and $G$-invariant, the set $M\setminus
A(f)$ is also $G$-invariant, and therefore we can consider the
following exact sequence of a pair:
$$...\rightarrow H^l(M/G,(M\setminus
A(f))/G)\stackrel{\alpha}{\rightarrow}
H^l(M/G)\stackrel{\beta}{\rightarrow} H^l((M\setminus
A(f))/G)\rightarrow....$$ By Lemma \ref{Lemma} the vector bundle
$\xi_M$ has a non-vanishing section over $M\setminus A(f)$. Thus
$\beta(e_p(\xi_M))=0$. Therefore there exists a non-trivial
element $$\mu\in H^{m(p^k-1)}(M/G,(M\setminus A(f))/G)$$ such that
$\alpha(\mu)=e_p(\xi_M)$. Since we are working over coefficients
in a field $\mathbb Z_p$ there exists a corresponding non-trivial
element $\tilde\mu\in H_{m(p^k-1)}(M/G,(M\setminus A(f))/G)$. Then
by Alexander duality we have $$H^{\dim M-m(p^k-1)}(A(f)/G;\mathbb
Z_p)\neq 0,$$ and thus $\dim_{\mathbb Z_p} A(f)/G\ge \dim
M-m(p^k-1)$ (see \cite {AD}). Since the group $G$ is finite it
easily follows that
$$\dim A(f)\ge \dim_{\mathbb Z_p} A(f)\ge \dim M-m(p^k-1),$$ and
we are done.
\end{proof}

\section*{Acknowledgements}
I would like thank A.N. Dranishnikov and Y.B. Rudyak for valuable
discussions and important remarks.

\end{document}